\documentclass[11pt,reqno]{amsart}

\usepackage[T1]{fontenc}
\usepackage[utf8]{inputenc}
\usepackage{lmodern}
\usepackage{microtype}

\usepackage[letterpaper,margin=1.1in]{geometry}

\usepackage{amssymb}
\usepackage{mathtools}
\usepackage{nicefrac}
\usepackage{tabularx}
\usepackage{enumitem}
\usepackage{graphicx}
\usepackage{booktabs}

\usepackage[hidelinks]{hyperref}
\usepackage[nameinlink,noabbrev]{cleveref}
\hypersetup{
  pdftitle={Artin-Schreier Theory and Local Brauer Invariants in Characteristic Two},
  pdfauthor={Hasitha Geekiyanage},
  pdfsubject={Mathematics research article}
}
\theoremstyle{plain}
\newtheorem{theorem}{Theorem}[section]
\newtheorem{lemma}[theorem]{Lemma}
\newtheorem{proposition}[theorem]{Proposition}
\newtheorem{corollary}[theorem]{Corollary}

\theoremstyle{definition}
\newtheorem{definition}[theorem]{Definition}
\newtheorem{example}[theorem]{Example}

\theoremstyle{remark}
\newtheorem{remark}[theorem]{Remark}

\numberwithin{equation}{section}
\crefname{assumption}{assumption}{assumptions}
\Crefname{assumption}{Assumption}{Assumptions}
\crefname{question}{question}{questions}
\Crefname{question}{Question}{Questions}

\newcommand{\Z}{\mathbb{Z}}
\newcommand{\Q}{\mathbb{Q}}

\newcommand{\F}{\mathbb{F}}
\newcommand{\G}{\mathfrak{G}}

\DeclareMathOperator{\ch}{char}
\DeclareMathOperator{\sep}{sep}
\DeclareMathOperator{\Hom}{Hom}

\DeclareMathOperator{\Gal}{Gal}
\DeclareMathOperator{\Br}{Br}
\DeclareMathOperator{\OP}{op}
\DeclareMathOperator{\im}{im}
\DeclareMathOperator{\Tr}{Tr}
\DeclareMathOperator{\N}{N}
\DeclareMathOperator{\Nrd}{Nrd}
\DeclareMathOperator{\Endo}{End}
\DeclareMathOperator{\inv}{inv}

\title[Artin-Schreier Theory]{Artin--Schreier Theory and Local Brauer Invariants in Characteristic Two}
\author[]{Hasitha Geekiyanage}
\address{Department of Mathematics and Statistics,
Texas Tech University, Lubbock, Texas, USA}
\date{}

\begin{document}

\begin{abstract}
These are survey notes on Artin-Schreier theory with a focus on separable quadratic extensions in characteristic two. Starting with some Galois theory, we classify separable quadratic extensions, then compute their norm forms, and construct the associated cyclic algebras; such an algebra has zero Brauer class exactly when its parameter is a norm, or equivalently, when the associated conic has a rational point. Finally, we consider quadratic constant extensions of global function fields and describe how the degree of a place and the parity of a valuation determine local splitting and the local Brauer invariant.
\end{abstract}

\maketitle

\tableofcontents

% ----------------------------------------------------------------------
\section{Introduction}
The setup is a field $K$ of characteristic $p>0$. The polynomial $X^n-a$ over $K$ used in classical Kummer theory is purely inseparable and hence fails for Galois theory when $n=p$. In this case, the Artin-Schreier polynomial $X^p-X-a$ is used. In particular, in characteristic two, a separable quadratic extension is generated by a root $\alpha$ satisfying $\alpha^2+\alpha=a$.

The main references are \cite{Lang}, \cite{Serre}, \cite{Voight} and \cite{Milne}.

Chapter VI, \S\S4-6 of \cite{Lang} used as the reference for trace, Hilbert 90, and cyclic extensions while \S8 used as for Artin-Schreier pairing for abelian extensions of exponent $p$.

Galois cohomology referred to chapter X of \cite{Serre}; ramification and norms to chapters IV-V; and Brauer group of a local field to chapter XII.

Almost all of quaternion algebra facts used follows from \cite{Voight} whereas Brauer classes, local and global fields are closely related to \S8.3 and chapters 13-14 respectively of the same text.

Some results from local class field theory and cohomological description of Brauer groups used are referred to chapters III-IV of \cite{Milne}.

In addition, some standard results appears in \cite{Poonen} are cited without proofs in later sections. \\
\textit{Remark}: This paper is self-contained and all the results presented are known, up to specialization to particular cases.

\textbf{Conventions:} All fields have characteristic $p>0$ where $p$ is a prime number. Field homomorphisms preserve identities as usual. Characteristic of the field $K$ is sometimes denoted as $\ch K$; and $K^{\sep}$ denotes the separable closure contained in the algebraic closure $\overline{K}$. The center of the algebra $A$ is denoted by $Z(A)$.

% ----------------------------------------------------------------------
\section{The Artin–Schreier map and cyclic extensions}

Let $K$ be a field with $\ch K= p>0$.

\begin{definition}
    Define the \emph{Artin–Schreier map} by $\wp: K\longrightarrow K$ by $z\longmapsto z^p-z$. The image of this map $\wp(K)=\{z^p-z: z\in K\}$ is called the Artin–Schreier group of $K$. 
\end{definition}

The kernel of this map is $\mathbb{F}_p$. The quotient $K/\wp(K)$ is an $\mathbb{F}_p$-vector space. If $\ch K=2$ then Artin–Schreier map becomes $\wp(z)=z^2+z$.

\begin{example}[$\wp$  is not multiplicative in general]
    Let $K=\mathbb{F}_2(t)$. Set $x=1,y=t$. Then $\wp(xy)=\wp(t)=t^2+t\neq0$. But $\wp(x)\wp(y)=(1^2+1)(t^2+t)=0$.
\end{example}

\begin{example}[$\wp(K)\subseteq K$ is not an ideal]
    Let $K=\mathbb{F}_4=\mathbb{F}_2(\omega)$ where $\omega^2+\omega+1=0$. Then $\wp(0)=\wp(1)=0$ and $\wp(\omega)=\wp(\omega+1)=1$. So, $\wp(K)=\{0,1\}\subseteq K$ is not an ideal.
\end{example}

\begin{theorem}\label{Artin–Schreier polynomial}
    Given $a\in K$, the polynomial $f_a(X)=X^p-X-a$ is separable. Moreover, one and only one of the following holds:
    \begin{enumerate}
        \item $a\in \wp(K)$, and $f_a$ splits completely over $K$.
        \item $a\notin \wp(K)$, and $f_a$ is irreducible; if $\alpha$ is a root then $K(\alpha)$ is cyclic over $K$ of degree $p$ with automorphism $\alpha\mapsto\alpha+c$ for each $c\in \mathbb{F}_p$.
    \end{enumerate}
\end{theorem}
This is the part (ii) of the Theorem 6.4 in \cite[Chap.~VI, \S6]{Lang}.
\begin{proof}
    $f_a$ is separable as its derivative is $-1\neq 0$. For any root $\alpha\in K^{\sep}$, roots of $f_a$ are precisely $\alpha+c$, $c\in\mathbb{F}_p$ as $(\alpha+c)^p-(\alpha+c)-a=\alpha^p-\alpha-a=0$. Hence $K(\alpha)$ contains all the roots and is a finite Galois extension of $K$. The map $\Gal(K(\alpha)/K)\longrightarrow(\mathbb{F}_p,+)$ given by $\sigma\longmapsto\sigma(\alpha)-\alpha$ is an injective group homomorphism. Since the group $\mathbb{F}_p$ has prime order, image of this map is either $0$ or $\mathbb{F}_p$. Thus, $[K(\alpha):K]$ is either $1$ or $p$.

    In the former case $\alpha\in K$, so $a\in\wp(K)$ and all roots are in $K$. In the latter case the minimal polynomial of $\alpha$ is $f_a$ and $\Gal(K(\alpha)/K)\cong(\mathbb{F}_p,+)$.
\end{proof}

The following additive form of Hilbert's 90 theorem provides a connection between cyclic extensions and their corresponding trace maps. Also see Theorem 6.3 in \cite[Chap.~VI, \S6]{Lang}

\begin{lemma}[Hilbert's 90 theorem, additive form]\label{additive Hilbert 90}
    Let $L/K$ be a finite cyclic Galois extension. Suppose $\Gal (L/K)=\langle\sigma\rangle$. Then $\ker (\Tr_{L/K})=\{\sigma(z)-z:z\in L\}$.
\end{lemma}
\begin{proof}
    Let $n=[L:K]$. The $K$-linear map $\Delta_\sigma:L\longrightarrow L$, $\Delta_\sigma(z)=\sigma(z)-z$ has kernel $K$. So, $\dim_K \im (\Delta_\sigma)=n-1$ and 
    $$\Tr_{L/K}(\Delta_\sigma(z))=\sum_{r=0}^{n-1}(\sigma^{r+1}(z)-\sigma^r(z))=0.$$
    Conversely, the trace map of a finite separable extension is non-zero \cite[Chap.~VI, \S5, Theorem~5.6]{Lang}. So, $\Tr_{L/K}:L\longrightarrow K$ is surjective and hence its kernel has dimension $n-1$. Therefore, $\im(\Delta_\sigma)=\ker\Tr_{L/K}$.
\end{proof}

\begin{theorem}[Artin–Schreier]\label{Artin–Schreier theorem}
    A field extension $L/K$ is cyclic of degree $p=\ch K$ if and only if $L=K(\alpha)$ where $\alpha^p-\alpha=a$ for some $a\in K\setminus\wp(K) $.
\end{theorem}
\begin{proof}
    Suppose $\Gal(L/K)=\langle\sigma\rangle$ has order $p$. Since $\Tr_{L/K}(1)=p=0$ by Lemma \ref{additive Hilbert 90} there is $\alpha\in L$ such that $\sigma(\alpha)-\alpha=1$. Now, $\sigma(\alpha^p-\alpha)=(\alpha+1)^p-(\alpha+1)=\alpha^p-\alpha$. So, $a:=\alpha^p-\alpha\in K$. The elements $\alpha,\alpha+1,\dots,\alpha+(p-1)$ are distinct conjuhates and hence $[K(\alpha):K]=p$ and $L=K(\alpha)$. If $a\in\wp(K)$ then by Theorem \ref{Artin–Schreier polynomial} we get $\alpha\in K$ which is a contradiction. 

    The reverse implication follows from Theorem \ref{Artin–Schreier polynomial}.
\end{proof}

\begin{remark}
    When $p>2$, degree $p$ separable extensions need not be Galois. When $p=2$, separable quadratic extensions are galois and hence cyclic, so the Theorem \ref{Artin–Schreier theorem} covers such extensions. 
\end{remark}

% ----------------------------------------------------------------------
\section{Parameterization of Artin–Schreier Extensions} 

For $a\in K$, define the $K$-algebra
$$E_a:=\nicefrac{K[X]}{(f_a(X))}$$
where $f_a(X)=X^p-X-a$. Since the derivative of $f$ is non-zero, it has distict roots in $K^{\sep}$. Thus,
$$E_a\otimes_KK^{\sep}\cong \nicefrac{K^{\sep}[X]}{(f_a(X))}\cong\prod_{c\in \mathbb{F}_p}K^{\sep}.$$
That is, $E_a$ is an  étale $K$-algebra. More explicitly,

\[
 E_a\cong
 \begin{cases}
 K^p;& \text{if}\hskip1mm a\in\wp(K)\\
 \text{a cyclic field extension of degree }p;& \text{if} \hskip1mma\notin\wp(K).
 \end{cases}
\]
where $K^p$ denotes the direct sum of $p$ copies of $K$.

Note that the parameter $a$ in $\alpha^p-\alpha=a$ is not uniquely determined by the extension $K(\alpha)$. Indeed, if $d\in K$ and $\beta=\alpha+d$ then $\beta^p-\beta=(\alpha+d)^p-(\alpha+d)=a+\wp(d)$; similarly, if $c\in\mathbb{F}_p^\times$ and $\beta=c\alpha$ then $\beta^p-\beta=(c\alpha)^p-c\alpha=c^p\alpha^p-c\alpha=c(\alpha^p-\alpha)=ca$. The following proposition states that a change of parameter $a\mapsto b$ gives the same non-trivial Artin-Schreier extension exactly when $b=ca+\wp(d)$ for some $c\in\mathbb{F}_p^\times$ and $d\in K$.

\begin{proposition}
    Let $a,b\notin\wp(K)$ and let $L_a=K(\alpha),L_b=K(\beta)\subseteq K^{\sep}$ where $\wp(\alpha)=a$ and $\wp(\beta)=b$. Then $L_a=L_b$ if and only if $b=ca+\wp(d)$ for some $c\in\mathbb{F}_p^\times$ and $d\in K$. Therefore, there is a one-to-one correspondence between degree $p$ cyclic subextensions of $\nicefrac{K^{\sep}}{K}$ and one dimensional $\mathbb{F}_p$-subspaces of $K/\wp(K)$.
\end{proposition}
This is a special case of Theorem 8.3 in \cite[Chap.~VI, \S8]{Lang} considering cyclic extensions of degree $p$.
\begin{proof}
    Suppose $L_a=L_b=L$. Choose $\sigma\in\Gal(L/K)$ with $\sigma(\alpha)=\alpha+1$. Since $\sigma$ fixes $\beta^p-\beta\in K$ and $X^p-X$ has all roots in $\mathbb{F}_p$, write $\sigma(\beta)-\beta=c$ for some non-zero $c\in\mathbb{F}_p$. This yields $\sigma(\beta-c\alpha)=\beta-c\alpha$ giving $d:=\beta-c\alpha\in K$. Applying $\wp$ we get $b=ca+\wp(d)$. 

    Conversely, if $b=ca+\wp(d)$ then $c\alpha+d$ is also a root of $f_b=X^p-X-b$. So, $\beta\in K(\alpha)$.
\end{proof}

Since $\mathbb{F}_2^\times=\{1\}$, each one dimensional $\mathbb{F}_2$-subspace has exactly one non-zero element. If $L_1,L_2\subseteq K^{\sep}$ are $K$-isomorphic separable quadratic extensions, then $L_1=L_2$ as both extensions are normal over $K$ and hence any $K$-embedding of $L_1$ into $K^{\sep}$ has image $L_1$. This yield the following corollary. 

\begin{corollary}
    If $\ch K=2$ then there is a one-to-one correspondence between the non-zero elements of $K/\wp(K)$ and separable quadratic extensions of $K$, up to $K$-isomorphism. The zero equivalence class corresponds to the étale algebra $K\times K$. That is, if $a\in\wp(K)$ then $X^2+X+a$ splits into two linear factors over $K$ and $\nicefrac{K[X]}{(X^2+X+a)}\cong K\times K$.  
\end{corollary}

\begin{lemma}[pole-order test]\label{lemma:1}
    Let $v:K^\times\longrightarrow\mathbb{Z}$ be a discrete valuation. If $v(a)<0$ and $p\nmid v(a)$, then $a\notin\wp(K)$.
\end{lemma}
\begin{proof}
    Suppose $a=z^p-z$ for some $z\in K$. If $v(z)\geq0$ then $v(a)=v(z^p-z)\geq\min\{pv(z),v(z)\}\geq0$. If $v(z)<0$ then $v(a)=v(z^p-z)=\min\{pv(z),v(z)\}=pv(z)$, a contradiction. 
\end{proof}

\begin{example}
    In $K=\mathbb{F}_p(t)$, the valuation at infinity gives $v_\infty(t)=-1<0$; but $p\nmid v_\infty(t)$. So, $t\notin\wp(K)$ and hence $X^p-X-t$ is irreducible over $K$. When $p=3$, the parameters $t,2t\notin\wp(K)$ have distinct non-zero classes in modulo $\wp(K)$ as $t-2t=-t=2t\notin\wp(K)$; but they give the same field because $\alpha^3-\alpha=t$ implies $(2\alpha)^3-2\alpha=2t$.
\end{example}

\begin{theorem}[finite Artin-Schreier pairing]
Let $V\subseteq K/\wp(K)$ be an $r$-dimensional $\F_p$-subspace, with basis $a_1,\ldots,a_r$. Choose $\alpha_i^p-\alpha_i=a_i$, and let $ L=K(\alpha_1,\ldots,\alpha_r)$. Then $L/K$ is Galois and
$$\Gal(L/K)\cong\Hom_{\F_p}(V,\F_p);\qquad [L:K]=p^r.$$
The pairing is given by $(\sigma,\overline a)\longmapsto\sigma(\alpha)-\alpha$, where $\wp(\alpha)=a$.
\end{theorem}
This is the finite-dimensional case of \cite[Chap.~VI, \S8, Theorem~8.3]{Lang}.

\begin{proof}
The compositum of the splitting fields of the $P_{a_i}$ is Galois. Translation of the chosen roots gives an injection
$$ G:=\Gal(L/K)\hookrightarrow\F_p^r.$$
If the image is a proper subspace then choose a nonzero vector $(c_1,\ldots,c_r)\in\F_p^r$ such that $\sum_i c_i(\sigma(\alpha_i)-\alpha_i)=0$
for every $\sigma\in G$. Then $\beta=\sum_i c_i\alpha_i$ would lie in $K$ and
$$\sum_i c_i a_i=\beta^p-\beta\in\wp(K).$$
This contradicts the fact that classes $\overline a_i$ are independent. Thus $G=\F_p^r$. Changing $\alpha$ by an element of $\F_p$, or changing the representative $a$ by $\wp(d)$ with $d\in K$, does not change the pairing. This completes the proof.
\end{proof}

% ----------------------------------------------------------------------
\section{A Cohomological description}
This section relates to chapter X of \cite{Serre}. We will discuss Brauer group again in section \ref{Brauer groups}. Also refer \cite[Chap.~X, \S\S1, 3--5]{Serre} and \cite[Chap.~IV, \S3]{Milne} for Galois cohomology and Brauer groups respectively.

Let $\G_K=\Gal (K^{\sep}/K)$ be the absolute Galois group. The Artin-Schreier sequence
\begin{equation}\label{Artin-Schreier sequence}
 0\longrightarrow\F_p\longrightarrow K^{\sep}
 \xrightarrow{\wp}K^{\sep}\longrightarrow0
\end{equation}
is exact additive sequence of discrete $\G_K$-modules. The map $\wp$ is surjective because each polynomial $X^p-X-a$ over $K^{\sep}$ has a root there. The additive form of Hilbert 90 in Galois cohomology states that $H^1(\G_K,K^{\sep})=0$. The long exact sequence corresponding to (\ref{Artin-Schreier sequence}) yields
\begin{equation}\label{cohomology eq}
K/\wp(K)\cong H^1(\G_K,\F_p)=\Hom_{\mathrm{cont}}(\G_K,\F_p)
\end{equation}
Here, $\G_K$ acts trivially on $\F_p$; so a 1-cocycle is a continuous homomorphism, and all 1-coboundaries are zero. This isomorphism can be described explicitly. Given $a\in K$, choose $\alpha\in K^{\sep}$ with $\wp(\alpha)=a$, and set
$$\chi_a(\sigma)=\sigma(\alpha)-\alpha \in\F_p$$
Notice that the character $\chi_a(\sigma)$ does not depend on the root chosen, because any two roots differ by an element of $\F_p$, which $\G_K$ fixes. Also $\chi$ will not change when $a$ is replaced by $a+\wp(d)$, since the root can be replaced by $\alpha+d$.

The isomorphism (\ref{cohomology eq}) can be verified as follows: the character $\chi$ is zero precisely when $\alpha$ is fixed by $\G_K$, equivalently when $a\in\wp(K)$. This proves the injectivity on the quotient. Now, every non-zero continuous character $\G_K\longrightarrow\F_p$ gives a cyclic extension of degree $p$, namely the fixed field of its kernel. Chose a generator on which the charater takes the value $1$. The construction in the proof of theorem \ref{Artin–Schreier theorem} gives an Artin-Schreier root. The zero character corresponds to $a=0$, so the map is surjective. 

The Brauer group has the following standard cohomological description. 
$$\Br(K)\cong H^2(\G_K,(K^{\sep})^\times)$$
For the proof see \cite[Chap.~IV]{Milne}. We will construct Brauer classes using central simple algebras in section \ref{Brauer classes}. Cohomology is not required to prove their splitting criterion.

Over a characteristic $p$ field, we have $X^p-1=(X-1)^p$; which shows that $1$ is the only $p$-th root of unity. Thus, classical separable Kummer theory cannot be used to describe cyclic extensions of degree $p$. Artin-Schreier theory instead theory instead uses the additive group $\F_p$. So, classical Kummer theory does not replace the sequence (\ref{Artin-Schreier sequence}). For basic Kummer theory see \cite[Chap.~IV, \S3]{Neukirch}. 

% ----------------------------------------------------------------------
\section{Artin–Schreier equations over Finite fields}
Let $q=p^r$ be a prime power, For $a\in \mathbb{F}_q$, the absolute trace for the extension $\mathbb{F}_q/\mathbb{F}_p$ is given by 
$$\Tr_{\mathbb{F}_q/\mathbb{F}_p}(a)=a+a^p+\cdots+a^{p^{r-1}}\in \mathbb{F}_p$$
this can be used to determine exactly which Artin-Schreier equations have solutions in $\mathbb{F}_q$. See \cite[Chap.~VI, \S5]{Lang} for more details about finite fields; see  \cite[\S12.3.11]{Voight} for the irreducibility in characteristic two case of the following proposition.

\begin{proposition}\label{kernel of trace}
    We have $\wp(\mathbb{F}_q)=\ker\Tr_{\mathbb{F}_q/\mathbb{F}_p}$ and $\nicefrac{\mathbb{F}_q}{\wp(\mathbb{F}_q)}\cong\mathbb{F}_p$. Moreover, $X^p-X-a$ is irreducible over $\mathbb{F}_q$ if and only if $\Tr_{\mathbb{F}_q/\mathbb{F}_p}(a)\neq0$.
\end{proposition}

\begin{proof}
    First note that every element in the image has zero trace as for any $z\in \mathbb{F}_q$, $\Tr(z^p-z)=\sum_{i=0}^{r-1}(z^p-z)^{p^i}=z^{p^r}-z=0$. 

    The kernel of $\wp$ has $p$ elements, so its image has $q/p$ elements. The trace is a non-zero $\mathbb{F}_p$-linear map ans defines a polynomial with degree $p^{r-1}<q$, so it cannot vanish at all $q$ elements of $\mathbb{F}_q$. Since the absolute trace map is surjective, its kernel has dimension $r-1$ and hence $q/p$ elements. Since the sets are finite this forces $\wp(\mathbb{F}_q)=\ker\Tr_{\mathbb{F}_q/\mathbb{F}_p}$. Again surjectivity of absolute trace map gives $\nicefrac{\mathbb{F}_q}{\wp(\mathbb{F}_q)}\cong\mathbb{F}_p$. Irreducibility condition folllows from Theorem \ref{Artin–Schreier polynomial}.
\end{proof}

\begin{example}
    The polynomial $X^2+X+1$ is irreducible over $\mathbb{F}_2$ and adjoining a root gives $\mathbb{F}_4$. But it splits over $\mathbb{F}_4$ as $1=\wp(\omega)$. In general, $\Tr_{\mathbb{F}_{2^r}/\mathbb{F}_2}(1)$ is $1$ if $r$ is odd and $0$ if $r$ is even. Thus, the class of $1$ is non-zero in $\nicefrac{\mathbb{F}_{2^r}}{\wp(\mathbb{F}_{2^r})}$ exactly when $r$ is odd. So, adjoining a root of $X^2+X+1$ does not give a quadratic extension of every finite field of characteristic $2$.
\end{example}

% ----------------------------------------------------------------------
\section{Characteristic two fields and Quadratic extensions}
From this section onward we assume $\ch K=2$ unless stated otherwise.

\subsection{Separable extensions, trace and norm equations}

Consider a separable quadratic extension $L/K$ with $L=K(\theta)$. The minimal polynomial $f$ of $\theta$ has the form $f(X)=X^2+cX+d$ where $c\neq 0$. Since $f(\theta)=0$, dividing this equation by $c^2\neq0$ and setting $\alpha=\theta/c$ we get $\alpha^2+\alpha=a$ where $a=d/c^2$. Non-identity automorphism for this extension is $\alpha\mapsto\alpha+1$. This gives Artin-Schreier description for a separable quadratic extension of characteristic 2. Note that this is the quadratic case of Theorem \ref{Artin–Schreier theorem}. For a non-trivial extension one must have $a\notin \wp(K)$. For more details on separable quadratic extensions over characteristic two fields see \cite[\S6.1]{Voight}.

In contrast, if $d\notin K^2$, then $X^2-d$ is irreducible over $K$ and has only one root over $\overline{K}$ with multiplicity two. Thus $K(\sqrt{d})/K$ is a purely inseparable of degree two.

\begin{proposition}\label{norm form}
    Let $a\notin \wp(K)$ and write $L=K(\alpha)$ where $\alpha^2+\alpha=a$. For $x,y\in K$
    $$\Tr_{L/K}(x+y\alpha)=y \hskip2mm\text{and}\hskip2mm \N_{L/K}(x+y\alpha)=x^2+xy+ay^2.$$ 
    This norm formula holds in the quadratic étale algebra $E_a=\nicefrac{K[X]}{(f(X))}$ where $f(X)=X^2+X+a$ regardless of the irreducibility of $f$.
\end{proposition}
\begin{proof}
    The non-identity automorphism maps $x+y\alpha$ to $x+y(\alpha+1)$; taking the sum of these two conjugates we get the trace. For the norm formula, note that $(x+y\alpha)(x+y(\alpha+1))=x^2+xy+y^2(\alpha^2+\alpha)=x^2+xy+ay^2$. 
    
    Since this uses only the relation $\alpha^2+\alpha=a$, the same calculation works in $E_a$. If $f$ splits then the norm is the multiplication of the two coordinates in $K\times K$.
\end{proof}

\begin{example}\label{example1}
    In the extension $\nicefrac{\mathbb{F}_4(t)}{\mathbb{F}_2(t)}$ we have $N(t+\omega)=(t+\omega)(t+\omega+1)=t^2+t+1$. So, $x^2+xy+y^2=t^2+t+1$ has a solution $(x,y)=(t,1)$.
\end{example}

Next, we give the multiplicative form of Hilbert's 90 theorem which identifies the elements of norm one. More general versions of this theorem can be found in \cite[Chap.~VI, \S6, Theorem~6.1]{Lang} and \cite[Chap.~X, \S1]{Serre}.

\begin{proposition}[Hilbert's 90 theorem, multiplicative form]
    Suppose $L/K$ is a separable quadratic extension with non-trivial automorphism $\sigma$. Then for any $u\in L$, $\N_{L/K}(u)=1$ if and only if $u=\frac{v}{\sigma(v)}$ for some $v\in L^\times$.
\end{proposition}
\begin{proof}
   Suppose $\N_{L/K}(u)=u\sigma(u)=1$. Now, $\sigma(u)=u^{-1}$. If $u=1$, set $v=1$; so $\frac{v}{\sigma(v)}=1=u$. If $u\neq1$, set $v=1+u\neq0$; so $\frac{v}{\sigma(v)}=\frac{1+u}{1+\sigma(u)}=\frac{1+u}{1+u^{-1}}=u$.
   
   Conversely, if $u=\frac{v}{\sigma(v)}$ then $\N_{L/K}(u)=u\sigma(u)=\left(\frac{v}{\sigma(v)}\right)\left(\frac{\sigma(v)}{\sigma^2(v)}\right)=1$.
\end{proof}

\subsection{Rational points on Projective conics} For $X,Y\in K$ by the norm formula $\N_{L/K}(X+Y\alpha)=X^2+XY+aY^2$. Thus, for $a\in K$ and $b\in K^\times$, the projective conic $\mathrm{C}_{a,b}$ in $\mathbb{P}^2_k$ is given by the equation $X^2+XY+aY^2=bZ^2$. A $K$-rational point is a point $[X:Y:Z]$ in $\mathbb{P}^2_k$ which satisfies this equation.

\begin{proposition}\label{projective conic}
    The projective conic $\mathrm{C}_{a,b}$ is smooth; and the set $\mathrm{C}_{a,b}(K)$ of $K$-rational points on it is non-empty if and only if $b\in\N_{E_a/K}(E^\times_a)$. 
\end{proposition}
\begin{proof}
    In characteristic 2, the homogeneous polynomial $F=X^2+XY+aY^2-bZ^2$ has partial derivatives $F_X=Y, F_Y=X$ and $F_Z=0$. For a singular point we must have $X=Y=0$; but this forces $Z=0$ which is not a projective point.

    If $b=\N(x+y\alpha)$ then $[x:y:1]$ is a rational point. Conversely suppose $\mathrm{C}_{a,b}(K)\neq\varnothing$. \\
    Case I ($E_a$ is a field): 
    Let $[X:Y:Z]\in\mathrm{C}_{a,b}(K)$. If $Z=0$ then $\N_{E_a/K}(X+Y\alpha)=X^2+XY+aY^2=0$. Because the projective coordinates cannot be zero simultaneously, we get $X+Y\alpha\neq0$. But this is impossible as non-zero element cannot have zero norm. So, $Z\neq0$. Hence $[x:y:1]\in\mathrm{C}_{a,b}(K)$ where $x=\frac{X}{Z}$ and $y=\frac{Y}{Z}$. Consequently, $b=x^2+xy+ay^2=\N_{E_a/K}(x+y\alpha)$ with $x+y\alpha\in E_a^\times$. \\
    Case II ($E_a\cong K\times K$): In this case norm map given by $(u,v)\longmapsto uv$ is surjective and hence any $b\in K^\times$ is a norm; by previous implication the conic has a rational point. 
\end{proof}

For general correspondences between conics and quaternion algebras see \cite[Corollary~6.4.8 and Theorem~6.4.11]{Voight}. Smooth projective conic is a one-dimensional \emph{Severi-Brauer} variety; \cite[Chap.~X, \S6]{Serre} defines these varieties and explains their algebraic representation via central simple algebras. 

% ----------------------------------------------------------------------
\section{Quaternion algebras and Brauer groups}\label{Brauer classes}

\subsection{Quaternion algebras of characteristic two}

\begin{definition}
    Let $a\in K$ and $b\in K^\times$. The quaternion algebra $Q(a,b)$ associated to $a$ and $b$ is the 4-dimensional $K$-algebra with basis $\{1,i,j,ij\}$ subject to the relations
    $$i^2+i=a,\qquad j^2=b,\qquad ji=(i+1)j.$$
\end{definition}

Note that this replaces the familiar noncommutativity condition used in the case when characteristic not equal to two; as in characteristic two $ij+ji=0$ implies commutativity.

A finite dimensional $K$-algebra is called central simple if its center is $K$ and has no nontrivial proper ideals. 

\begin{proposition}
    The $K$-algebra $Q(a,b)$ is central simple. For any field extension
    $E/K$ containing a root of $X^2+X+a$, the $E$-algebra
    $Q(a,b)\otimes_K E$ is isomorphic to $M_2(E)$.
\end{proposition}

\begin{proof}
    Let $\alpha\in E$ be a root of $X^2+X+a$. Over $E$, let
    \[
    I=\begin{pmatrix}
        \alpha & 0\\
        0 & \alpha+1
    \end{pmatrix},
    \qquad
    J=\begin{pmatrix}
        0 & b\\
        1 & 0
    \end{pmatrix}.
    \]
    Then these satisfy $I^2+I=aI_2$, $J^2=bI_2$, and $JI=(I+I_2)J$,
    where $I_2$ is the identity matrix. So we get an $E$-algebra
    homomorphism $Q(a,b)\otimes_K E\longrightarrow M_2(E)$.

    It is easy to see that the matrices $I_2,I$ span the diagonal
    matrices in $M_2(E)$, while $J,IJ$ span the off-diagonal matrices.
    Since $b\neq 0$ and $(\alpha+1)-\alpha=1\neq 0$, these four matrices
    form a basis of $M_2(E)$. They are the images of the basis
    $1,i,j,ij$ of $Q(a,b)\otimes_K E$, so it follows that
    $Q(a,b)\otimes_K E\cong M_2(E)$.

    Since $M_2(E)$ is simple, under this scalar extension any nonzero
    two-sided ideal of $Q(a,b)$ becomes all of $M_2(E)$. Scalar
    extension preserves dimension, so the ideal has dimension four
    over $K$ and is therefore all of $Q(a,b)$.

    Under the same isomorphism, the scalar extension
    $Z(Q(a,b))\otimes_K E$ is contained in $Z(M_2(E))=E I_2$, which
    is one-dimensional over $E$. Therefore $Z(Q(a,b))$ has dimension
    at most one over $K$. Since it contains $K$, this implies
    $Z(Q(a,b))=K$.
\end{proof}

If $a\notin\wp(K)$, then the subalgebra $K(i)$ is the quadratic field $L=K(\alpha)$. Writing $i=\alpha$ we get the cyclic algebra described by
$$Q(a,b)=L\oplus Lj,\qquad jl=\sigma(l)j,\qquad j^2=b$$
where $\sigma(\alpha)=\alpha+1$ and $l\in L$.
In the split case, one may replace $L$ by quadratic étale algebra $E_a$ and the same formula works. 

\subsection{Brauer groups}\label{Brauer groups}
Two central simple algebras $A$ and $B$ are Brauer equivalent if $M_m(A)\cong M_n(B)$ as $K$-algebras for some $m,n\geq1$. Define the Brauer group $\Br(K)$ of $K$ as the set of these equivalece classes. The group operation is induced by the tensor product:
$$[A]+[B]=[A\otimes_KB]$$
Matrix algebras over $K$ gives the zero class and the inverse of $[A]$ is $[A^{\OP}]$.

Artin-Wedderburn theorem \cite[Theorem~2.1.3]{GilleSzamuely} for simple central algebras implies that $[A]=0$ if and only if $A\cong M_n(K)$ for some $n\geq1$. See \cite[Chap.~IV]{Milne} for more about Brauer groups.

We use the notation 
$$[a,b)_K:=[Q(a,b)]\in\Br(K)$$
to denote the Brauer class, reserving $Q(a,b)$ to denote the algebra itself. Thus, $[a,b)_K=0$ means $Q(a,b)\cong M_2(K)$; such an algebra is called \emph{split}. 

There is a natural isomorphism $Q(a,b)\longrightarrow Q(a,b)^{\OP}$, known as \emph{quaternion conjugation}, given on generators by $i\mapsto i+1$ and $j\mapsto j$ and reversing products. This yields $2[a,b)_K=0$; or equivalently, $[a,b)_K$ is in the $2$-torsion subgroup $\Br(K)[2]$. 

Now for any $x,y\in E_a$ the conjugate of $q=x+yj$ is given by $\bar{q}=\sigma(x)+yj$. It follows that the reduced norm of $x+yj$ is given by 
$$\Nrd_{Q/K}(x+yj)=q\bar{q}=\N_{E_a/K}(x)+b\N_{E_a/K}(y).$$
If we write $x=x_0+x_1\alpha$ and $y=y_0+y_1\alpha$ this norm formula becomes $x_0^2+x_0x_1+ax_1^2+b(y_0^2+y_0y_1+ay_1^2)$.

\begin{theorem}[The splitting criterion]\label{splitting criterion}
    Let $a\in K, b\in K^\times$ and $E_a=\nicefrac{K[X]}{(X^2+X+a)}$. Then the following are equivalent:
    \begin{enumerate}
        \item $[a,b)_K=0$ in $\Br(K)$;
        \item $Q(a,b)\cong M_2(K)$;
        \item $b\in N_{E_a/K}(E_a^\times)$;
        \item there exist $x,y\in K$ satisfying $x^2+xy+ay^2=b$;
        \item $\mathrm{C}_{a,b}(K)\neq\varnothing$.
    \end{enumerate}
\end{theorem}
This is a special version of \cite[Theorem~6.4.11]{Voight}.

\begin{proof}
    $(1)\Leftrightarrow(2)$ follows by the definition. $(3)\Leftrightarrow(4)\Leftrightarrow(5)$ follows by propositions \ref{projective conic}. and \ref{norm form}. 

    If $a\in\wp(K)$, then $E_a\cong K\times K$ and the norm is
    $\N(x,y)=xy$. So every $b\in K^\times$ is a norm. Also,
    $Q(a,b)\cong M_2(K)$ by Proposition~7.2, since $X^2+X+a$
    has a root in $K$. Thus both $(2)$ and $(3)$ hold in this case.

    We may therefore assume $a\notin\wp(K)$, so $E_a=L$ is a
    quadratic field. Let $\alpha$ be the class of $X$ in $L$ and
    let $\sigma$ be the nontrivial $K$-automorphism of $L$.
    Suppose $(3)$ holds. Then there is $c\in L^\times$ with
    $b=\N_{L/K}(c)=c\sigma(c)$. For each $l\in L$, let $M_l$
    denote multiplication by $l$ and define $J(v)=c\sigma(v)$
    on $L$. Then both $M_l$ and $J$ are $K$-linear,
    $J^2=b\,\mathrm{id}_L$, and $JM_l=M_{\sigma(l)}J$.
    Thus $i\mapsto M_\alpha$ and $j\mapsto J$ define a unital
    $K$-algebra homomorphism
    $\varphi:Q(a,b)\longrightarrow\Endo_K(L)$.
    Since $Q(a,b)$ is simple, $\varphi$ is injective. Now both
    of these algebras are $4$-dimensional over $K$, and hence
    $\varphi$ is an isomorphism. This gives $(2)$.

    Conversely, assume $Q(a,b)\cong\Endo_K(V)$ with $\dim_K V=2$.
    The subfield $L\subset Q(a,b)$ acts on $V$, making it a
    one-dimensional $L$-vector space. Let $J$ be the operator
    corresponding to $j$. For $v\neq 0$, we have $Jv=cv$ for some
    unique $c\in L$. Now $J$ is invertible, and hence $c\neq 0$.
    The quaternion relation $jl=\sigma(l)j$ gives
    \[
        bv=J^2v=J(cv)=\sigma(c)Jv=\sigma(c)cv.
    \]
    Since $v\neq 0$, we get $b=\N_{L/K}(c)$. This gives $(3)$.
\end{proof}

Observe that, extending to a field over which $E_a$ splits, the quaternion algebra also splits; but this does not guarantee the existence of $K$-rational points.

Replacing either $i$ by $i+d$ where $d\in K$, or $j$ by $cj$ where $c\in E_a^\times$ gives new generators satisfying the quaternion algebra relations with new parameters. This yields the following isomorphisms:
$$Q(a+\wp(d),b)\cong Q(a,b),\qquad\qquad Q(a,b\N(c))\cong Q(a,b)$$
In particular, multiplying $b$ by any non-zero square preserves the second isomorphism; in fact, if $s\in K^\times$ then $\sigma(s)=s$ and hence $\N_{E_a/K}(s)=s^2$.
For a fixed $E_a$, the reverse implication can be obtained using the Skolem–Noether theorem, given in \cite[Corollary~7.7.6]{Voight}.

\begin{example}
    Since $1$ and every non-zero square are norms $[a,1)_K=0$ and $[a,c^2)_K=0$ for $c\in K^\times$. If $a\neq0$ then $\N(\alpha)=a$, so $[a,a)_K=0$. By example \ref{example1}, $[1,t^2+t+1)_{\mathbb{F}_2(t)}=0$.
\end{example}

Let $L/K$ be a field extension. We define the \emph{relative Brauer group} of this extension as
$$\Br(L/K):= \ker\!\left(\Br(K)\xrightarrow{\operatorname{res}_{L/K}}\Br(L)\right).$$
where $\operatorname{res}_{L/K}$ extends scalars. For $L=K(\alpha)$ with $\wp(\alpha)=a$ by cyclic algebra theorem (see \cite[\S~1.5.7]{Poonen}) we have the isomorphism
$$K^\times/\N_{L/K}(L^\times)\xrightarrow{\;\sim\;}\Br(L/K),\qquad b \longmapsto [a,b)_K.$$
From the norm criterion we proved, the kernel of the map $b\mapsto [a,b)_K$ is $N_{L/K}(L^\times)$. The fact that every element
of the relative Brauer group is represented by such a cyclic algebra follows from the subsection 1.5.7 of \cite{Poonen}. From Galois cohomology \cite[Chap.~IV, \S~3]{Milne}, this corresponds to the natural isomorphism 
$$ H^2\bigl(\Gal(L/K),L^\times\bigr)
\cong K^\times/\N_{L/K}(L^\times).$$
Above isomorphism gives the following identity of Brauer classes:
\[
[a,b_1b_2)_K=[a,b_1)_K+[a,b_2)_K.
\]

% ----------------------------------------------------------------------
\section{Artin-Schreier classes over Laurent series fields}

In this section, any prime characteristic $p$ is allowed. Let $q=p^r$ be a prime power. Take $K=\mathbb{F}_q((t))$. An element in $K$ is a Laurent series $f=\sum_{n\geq n_0}a_nt^n$ with only finitely many negative powers. For a non-zero $f$, its valuation $v(f)$ is the least exponent with non-zero coefficient. The valuation ring is $\mathcal{O}_K=\mathbb{F}_q[[t]]$, with maximal ideal $t\mathcal{O}_K$ and residue field $\mathbb{F}_q$. The field $K$ is complete for this valuation. 

\begin{lemma}\label{lemma 2}
    The restricted Artin-Schreier map $t\mathcal{O}_K\longrightarrow t\mathcal{O}_K$ is a bijection. For $h\in t\mathcal{O}_K$, a solution for $z^p-z=h$ is given by $z=-\sum_{n=0}^{\infty}h^{p^n}$.
\end{lemma}
\begin{proof}
    For any non-zero $h$ we have $v(h^{p^n})=p^nv(h)\longrightarrow\infty$. By completeness of $K$, above series converges. Note that, every term except $h$ cancels in $z^p-z$ proving surjectivity. Injectivity follows because $\ker\wp=\mathbb{F}_p$, and $\mathbb{F}_p\cap t\mathcal{O}_K=\{0\}$. 
\end{proof}

Choose $c\in\mathbb{F}_q$ with $\Tr_{\mathbb{F}_q/\mathbb{F}p}(c)=1$.

\begin{proposition}
    Each class in $\nicefrac{K}{\wp(K)}$ has a unique representative
    \begin{equation}\label{eq:2}
    \varepsilon c+\sum_{\substack{m>0\\p\nmid m}}c_m t^{-m},
    \qquad \varepsilon\in\F_p,\quad c_m\in\F_q,
    \end{equation}
    where only finitely many $c_m$ are non-zero.
\end{proposition}
\begin{proof}
    By lemma \ref{lemma 2}, write any Laurent series as $\sum_{m=1}^{M}a_{-m}t^{-m}+a_0+\wp(u)$ for some $u\in t\mathcal{O}_K$. Suppose $dt^{-pm}$ is a negative term. Since $\mathbb{F}_q$ is perfect there is $e\in\mathbb{F}_q$ with $e^p=d$. This gives $\wp(et^{-m})=dt^{-pm}-et^{-m}$. Repeat this operation always choosing the largest pole divisible by $p$; this process terminates. To reduce the constant term let $\varepsilon=\Tr(a_0)\in\mathbb{F}_p$. Since $\Tr(c)=1$ and trace is $\mathbb{F}_p$-linear $\Tr(a_0-\varepsilon c)=0$. By proposition \ref{kernel of trace}, $a_0\equiv\varepsilon c (\mathrm{mod} \hskip1mm\wp(K))$.

    To prove uniqueness, for two representatives belong to the same class, their difference is again has the same form. If it has a pole, its largest pole order is not divisible by $p$, by lemma \ref{lemma:1} it is not in $\wp(K)$. If the difference is a constant, say $\varepsilon c\in\wp(K)$. Taking the absolute trace gives $\varepsilon=0$.
\end{proof}

\begin{example}
    Let $\mathbb{F}_2((t))$. The expression $t^{-4}+t^{-2}=\wp(t^{-2})$ represents the zero class. Also, $t^{-6}+t^{-1}\equiv t^{-3}+t^{-1} (\mathrm{mod}\hskip1mm \wp(K))$. So, the reduced representative has the largest pole order $3$.
\end{example}

% ----------------------------------------------------------------------
\section{Unramified and totally ramified extensions}
Let $L/K$ be a finite extension of complete discretely valued fields, with normalized valuations $v$ on $K$ and $w$ on $L$. The \emph{ramification index} $e$ satisfies $w(x)=ev(x)$ for all $x\in K^\times$. the \emph{residue degree} $f$ is the degree of the extension of their residue fields. We say $L/K$ is \emph{unramified} if $e=1$, $f=[L:K]$ and the residue field extension is separable. We say $L/K$ is \emph{totally ramified} if $e=[L:K]$ and the two residue fields are same.

\begin{proposition}\label{reduced form proposition}
    Let $K=\mathbb{F}_q((t))$, and suppose $a\in K$ has the reduced form (\ref{eq:2}).
    \begin{enumerate}
        \item If $a=0$ then the  étale algebra $E_a$ splits.
        \item If $a=\varepsilon c$ with $\varepsilon\neq0$, adjoining a root gives the unramified extension $\nicefrac{\mathbb{F}_{q^p}((t))}{\mathbb{F}_q((t))}$.
        \item If $a$ has a pole of order $m>0$, then $p\nmid m$ and $L=K(\alpha)$ where $\alpha^p-\alpha=a$ is totally ramified of degree $p$; normalized valuation $w$ on $L$ satisfies $w(t)=p$ and $w(\alpha)=-m$.
    \end{enumerate}
\end{proposition}
\begin{proof}
    (1) follows from theorem \ref{Artin–Schreier polynomial}. For (2), since $c\in\mathbb{F}_q$ was chosen with absolute trace $1$, $\Tr_{\mathbb{F}_q/\mathbb{F}_p}(\varepsilon c)=\varepsilon\neq0$. So, $X^p-X-\varepsilon c$ is irreducible over $\mathbb{F}_q$. Adjoining a root gives $\mathbb{F}_{q^p}((t))$. Notice that, uniformizer remains $t$, so $e=1$, residue field extension gives $f=p$.

    For (3), pole-order test gives $a\notin\wp(K)$. By irreducibility criterion $[L:K]=p$ where $K=\mathbb{F}_q((t))$, $L=K(\alpha)$. Since $w|_K=ev$ we get $w(a)=-em<0$ implying $w(\alpha)<0$ and hence $pw(\alpha)=w(\alpha^p-\alpha)=w(a)=-em$. Since $p\nmid m$, we get $p\mid e$. But $e\leq [L:K]=p$. So, $e=p$. Now, $f=1$ and hence valuation identities are satisfied. 
\end{proof}

The pole order of the reduced parameter also gives the lower ramification break. The calculation below shows the relationship between Artin-Schreier equations and ramification groups discussed in \cite[Chap.~IV, \S1]{Serre}. The following proposition is a special case considered $K=\F_q((t))$ of the proposition 3.7.8 in \cite{Stichtenoth}.

\begin{proposition}
    In the case (3) of proposition \ref{reduced form proposition}, the unique lower ramification break is $m$. The relative ideal $\mathcal{D}_{L/K}$ has exponent $(p-1)(m+1)$ with respect to the normalized valuation $w$.
\end{proposition}
\begin{proof}
    Choose integers $u,s$ with $pu-ms=1$ and $1\leq s\leq p-1$. Then $\pi=t^u\alpha^s$ has valuation one; so it is a uniformizer of $L$. Any non-identity automorphism maps $\alpha$ to $\alpha+d$ where $d\in\F_p^\times$. Hence
    $$\sigma(\pi)-\pi=t^u\bigl((\alpha+d)^s-\alpha^s\bigr).$$
    The highest degree term of $\alpha$ on the right hand side is $sd\alpha^{s-1}$ whose coefficient is nonzero. Since $w(\alpha)<0$, this term has strictly smaller valuation than all the other terms. Therefore
    $$w(\sigma(\pi)-\pi)=pu-m(s-1)=m+1.$$
    For a totally ramified Galois extension, the lower group $G_i$ consists of automorphisms such that $w(\sigma(\pi)-\pi)\geq i+1$. Hence $G_0=\cdots=G_m=\Gal(L/K)$ and $G_{m+1}=1$. The relative formula $w(\mathfrak D_{L/K})=\sum_{i\geq0}(|G_i|-1)$, proved in \cite[Chap.~IV, \S1]{Serre}, gives \((p-1)(m+1)\).
\end{proof}
For $p=2$, we have the following classification:

\begin{center}
\small
\begin{tabularx}{\linewidth}{@{}XX@{}}
\toprule
Reduced Artin--Schreier parameter & Associated quadratic \'etale algebra \\
\midrule
Zero & \(K\times K\) \\
Nonzero constant class & Unramified quadratic field \\
At least one negative term of odd pole order & Totally ramified quadratic
field; the largest pole order is its lower break \\
\bottomrule
\end{tabularx}
\end{center}

% ----------------------------------------------------------------------
\section{Local Brauer invariants}

Let's come back to characteristic two. Consider the unramified quadratic extension $L/K$ where $K=\mathbb{F}_{q}((t))$, $q=2^r$ and $L=\mathbb{F}_{q^2}((t))$. Both fields have uniformizer $t$.

\begin{proposition}\label{norm subgroup}
    The norm subgroup is given by 
    $$\N_{L/K}(L^\times)=\{b\in K^\times:v(b)\in 2\mathbb{Z}\}$$
    In particular, every unit in $\mathcal{O}_K$ is a norm.
\end{proposition}
This is the quadratic case of \cite[\S13.2.5]{Voight}. We give a proof for the Laurent series fields.

\begin{proof}
    Let $x=t^nu\in L^\times$ with $u\in\mathcal{O}_L^\times$. Since $\N(t)=t^2$ we get $\N(x)=t^{2n}\N(u)$; so $v(N(x))=2n$ is even.
    Conversely, let $b\in K^\times$ have even valuation. Write $b=t^{2m}u$ with $u\in\mathcal{O}_K^\times$. Since $\N_{L/K}(t^m)=t^{2m}$, it is enough to show that $u$ is a norm. Write $\N=\N_{L/K}$.The finite-field norm
    \[
    \mathbb F_{q^2}^{\times}\longrightarrow\mathbb F_q^{\times},\qquad c\longmapsto c^{q+1},
    \]
    is surjective. Choose $z_1\in\mathbb F_{q^2}^{\times}$ such that
    $\N(z_1)\equiv u(\mathrm{mod}\hskip1mm t)$. We proceed by induction. 
    Suppose that $z_n\in\mathcal O_L^\times$ satisfies $\N(z_n)\equiv u(\mathrm{mod}\hskip1mm t^n)$. Then for some $d_n\in\mathbb{F}_q$,
    \[
    \frac{u}{N(z_n)}
    \equiv 1+d_nt^n(\mathrm{mod}\hskip1mm t^{n+1}).
    \]
    Since the residue-field trace is surjective, there is $c_n\in\mathbb{F}_{q^2}$ such that $\Tr_{\mathbb{F}_{q^2}/\mathbb{F}_q}(c_n)=d_n$. Set $z_{n+1}=z_n(1+c_nt^n)$. Now, $\N(1+c_nt^n)=(1+c_nt^n)(1+c_n^qt^n)\equiv 1+d_nt^n(\mathrm{mod}\hskip1mm t^{n+1})$. So, $\N(z_{n+1})\equiv u(\mathrm{mod}\hskip1mm t^{n+1})$ completing the induction. Now the differences $z_{n+1}-z_n$ have valuations at least $n$, so $(z_n)$ is Cauchy. Completeness of $L$ gives a limit $z\in\mathcal O_L^\times$. Continuity of the norm gives $\N(z)=u$. Therefore, $b=t^{2m}u=N(t^m z)$ completing the proof. 
\end{proof}

It is a standard fact that nonarchimedean local fields gives an isomorphism
$$\inv_K:\operatorname{Br}(K)\xrightarrow{\sim}\mathbb{Q}/\mathbb{Z}.$$
called local invariant map whose construction and proof studied in local class field theory; see \cite[Chap.~XII]{Serre} and \cite[Chap.~III, \S1; Chap.~IV, \S4]{Milne}. In particular, since a quaternion Brauer class has order dividing $2$, it has invariant either $0$ or $\frac{1}{2}$; and the invariant is zero precisely when it splits.

\begin{theorem}[Unramified parity formula]
    Consider the unramified quadratic extension $L=K(\alpha)$ where $\alpha^2+\alpha=a$. For every $b\in K^\times$ we have,
    $$\inv_K([a,b)_K)=\frac{v(b)}{2}\pmod{\mathbb{Z}}$$
    More generally, for $a\in\mathcal{O}_K$ denote its residue by $a_0\in\mathbb{F}_q$. Set $\tau=\Tr_{\mathbb{F}_q/\mathbb{F}_2}(a_0)$ which is either $0$ or $1$. Then
    \begin{equation}\label{parity formula}
    \inv_K([a,b)_K)=\frac{\tau\,v(b)}{2}\pmod{\mathbb{Z}}
    \end{equation}
\end{theorem}
\begin{proof}
    First suppose that $L/K$ is a nontrivial unramified quadratic extension. By proposition \ref{norm subgroup} and the splitting criterion $[a,b)_K=0\Leftrightarrow b\in \N_{L/K}(L^\times)\Leftrightarrow v(b)$ is even. A quaternion class has invariant $0$ when it splits and $\frac{1}{2}$ otherwise. Hence we get the first formula.
    
    Now since $a-a_0\in t\mathcal{O}_K$, lemma \ref{lemma 2}  gives $a-a_0\in\wp(K)$. So, $[a,b)_K=[a_0,b)_K$. If $\tau=0$ then by finite-field trace criterion, $a_0\in\wp(\mathbb{F}_q)$. Thus the quaternion algebra splits for every $b\in K^\times$, and its invariant is zero. If $\tau=1$, the polynomial $X^2+X+a_0$ gives the nontrivial quadratic extension of the residue field. Consequently, adjoining a root gives the nontrivial unramified quadratic extension of $K$, so the first formula applies. Considering both cases we have \ref{parity formula}.
\end{proof}

\begin{example}
Let \(K=\mathbb{F}_2((t))\). The parameter \(a=1\) gives
the unramified quadratic extension \(L=\mathbb{F}_4((t))\).
By parity formula we get:

\begin{center}
\begin{tabular}{@{}lcc@{}}
\toprule
Symbol & Local invariant & Reason \\
\midrule
\([1,t)_K\) & \(1/2\) & \(v(t)=1\) \\
\([1,t^2)_K\) & \(0\) & \(v(t^2)=2\) \\
\([1,1+t+t^2)_K\) & \(0\) & the parameter is a unit \\
\bottomrule
\end{tabular}
\end{center}

The second parameter in the last row can explicitly express as a norm: $N_{L/K}(1+t\omega)=1+t+t^2$ where $\omega^2+\omega+1=0$. However, over $K'=\mathbb{F}_4((t))$, the algebra $[1,t)_{K'}$ splits because $1=\wp(\omega)\in\wp(\mathbb{F}_4)$. This shows why the residue-field
trace factor in \ref{parity formula} is essential.
\end{example}

\begin{example}
For $K=\mathbb{F}_2((t))$ let $L=K(\alpha)$ where $\alpha^2+\alpha=t^{-1}$.
Then \(L/K\) is a totally ramified quadratic extension. However, 
$$\N_{L/K}(\alpha)=\alpha(\alpha+1)=t^{-1},\qquad [t^{-1},t^{-1})_K=0.$$
Thus $t^{-1}$ is a norm, but its valuation $v(t^{-1})=-1$ is odd. So, the valuation parity criterion does not hold for arbitrary ramified quadratic extensions.
\end{example}

% ----------------------------------------------------------------------
\section{Quadratic extensions of global function fields}

Let $F$ be a global function field of characteristic two with full constant field $\mathbb{F}_q$: i.e., $F$ has transcendence degree one over $\mathbb{F}_q$, is finite over a rational function field, and the elements of $F$ algebraic over $\mathbb{F}_q$ are precisely $\mathbb{F}_q$. Choose Choose $c\in\mathbb{F}_q$ with $\Tr_{\mathbb{F}_q/\mathbb{F}_2}(c)=1$. This gives a quadratic field extension
$$E=F\mathbb{F}_{q^2}=F(\alpha),\qquad \alpha^2+\alpha=c.$$
The assumption $\mathbb{F}_q$ is a full constant field ensures that $\mathbb{F}_{q^2}$ is not already contained in $F$.

For each place $v$ of $F$, let $\kappa(v)$ denote its finite residue field. The degree of $v$ is given by $d_v=[\kappa(v):\mathbb{F}_q]$. The completion $F_v$ has the same residue field $\kappa(v)$ and both these fields are of characteristic $2$. the following proposition connects finite field trace with local norm. For more about local and global theory see \cite[Chaps.~13--14]{Voight} and \cite[Chap.~XII]{Serre}.

\begin{proposition}
    At each place $v$,
    \[
    E\otimes_F F_v\cong
    \begin{cases}
    F_v\times F_v,&d_v\text{ even};\\
    \text{the unramified quadratic field over }F_v,&d_v\text{ odd}.\end{cases}
    \]
    Moreover, for each $b\in F^\times$
    \begin{equation}\label{eq:globalinv} \inv_v([c,b)_F)=\frac{d_v\,v(b)}2\pmod{\mathbb{Z}},
    \end{equation}
    Here, $\inv_v$ means restriction to $\Br(F_v)$ and then applying the local invariant map.
\end{proposition}
\begin{proof}
    Observe that, $\Tr_{\kappa(v)/\mathbb{F}_2}(c)=\Tr_{\mathbb{F}_q/\mathbb{F}_2}(\Tr_{\kappa(v)/\mathbb{F}_q}(c))=\Tr_{\mathbb{F}_q/\mathbb{F}_2}(d_vc)=d_v \pmod{2}$. By proposition \ref{kernel of trace}, $X^2+X+c$ splits when $d_v$ is even and is irreducible when $d_v$ is odd. Since its derivative is $1$, if $d_v$ is even Hensel's lemma lifts both roots to $F_v$. The resulting two distinct linear factors gives, by Chinese remainder theorem, $E\otimes_FF_v\cong F_v\times F_v$. If $d_v$ is odd adjoining a root over $F_v$ gives an unramified quadratic extension. Formula (\ref{parity formula}) gives (\ref{eq:globalinv}).
\end{proof}

Set $F=\mathbb{F}_2(t)$ and $c=1$. Each monic irreducible polynomial $P(t)$ defines a finite place of degree $\deg P$. For $0\neq b\in\mathbb{F}_2(t)$, valuation at $P$ is the exponent with which $P$ occurs in $b$. The place at infinity has degree one, with $v_\infty(f/g)=\deg g-\deg f$. For the class $[1,t)_F$, the invariants are as follows:
\begin{center}
\begin{tabular}{@{}lccc@{}}
\toprule
Place & Degree & \(v(t)\) & \(\inv_v([1,t)_F)\) \\
\midrule
\(t=0\) & \(1\) & \(1\) & \(1/2\) \\
Infinity & \(1\) & \(-1\) & \(1/2\) \\
All other places & \(\deg v\) & \(0\) & \(0\) \\
\bottomrule
\end{tabular}
\end{center}
At $t=0$, the quaternion algebra does not split. So, the conic $X^2+XY+Y^2=tZ^2$ has no $F$-rational points.

Put $b=t^2+t+1$; which is irreducible over $\mathbb{F}_2[t]$. It has only one finite zero corresponding to a degree-two place, where the constant extension splits. Also, $v_\infty(b)=-2$. So, all local invariants of $[1,b)_F$ vanish. We already have an explicit global norm equation $b=\N(t+\omega)$ and a rational point $[t:1:1]\in\operatorname{C}_{1,b}(F)$.

% ----------------------------------------------------------------------
\section{Reciprocity, local-global principle and Brauer--Manin obstructions} 

It is a standard result in global class field theory that for a global function field $F$, the following Brauer reciprocity sequence is exact:
\begin{equation}\label{reciprocity eq}
 0\longrightarrow\Br(F)\longrightarrow
 \bigoplus_v\Br(F_v)
 \xrightarrow{\ \sum_v\inv_v\ }\Q/\Z\longrightarrow0.
\end{equation}
See \cite[Remark~14.6.10 and equation~(14.6.11)]{Voight}; for global class field theory see \cite{ArtinTate}. For constant symbols, it is straightforward to see that,
$$\sum_v\inv_v([c,b)_F)=\frac{1}{2}\sum_vd_v\,v(b)=0 \pmod{\Z}$$
because a principal divisor has degree zero.

The class $[1,t)_F$ discussed above has two invariants equal to $1/2$; their sum is zero in $\Q/\Z$ as reciprocity requuires, although the algebra is not split. By injectivity in \ref{reciprocity eq}, a global class is split exactly when every local invariant vanishes.

Combining this injectivity with theorem \ref{splitting criterion} we get the local-global principle for these conics:
$$\operatorname{C}_{a,b}(F)\neq\varnothing\Longleftrightarrow\operatorname{C}_{a,b}(F_v)\neq\varnothing\hskip2mm\text{ for every } v.$$
When $b=t$ we already noted the absence of $F$-rational points. This is a local obstruction; thus this example agrees with local-global principle.

Let $X$ be a smooth projective geometrically integral variety over $F$. A fixed class $\mathcal{B}\in\Br(X)$ can be evaluated at local points $P_v\in X(F_v)$ defining the following paring:
$$((P_v),\mathcal{B})\longmapsto\sum_v\inv_v(\mathcal{B}(P_v))\in\Q/\Z$$
where only finitely terms are non-zero. let $\mathbb{A}_F$ denote the ring of adeles of $F$ which is the restricted product of all completions $F_v$ of $F$. The Brauer--Manin set consists of the adelic points orthogonal to every class:
$$X(\mathbb{A}_F)^{\Br}=\left\{(P_v):\sum_v\inv_v(\mathcal{B}(P_v))=0 \text{ for every } \mathcal{B}\in\Br(X)\right\}$$

By reciprocity $X(F)\subseteq X(\mathbb{A}_F)^{\Br}$. If $X(\mathbb{A}_F)\neq\varnothing$ but $X(\mathbb{A}_F)^{\Br}=\varnothing$, then there is a Brauer--Manin obstruction to rational points; for more details see \cite[\S8.2]{Poonen}.

Artin-Schreier symbols allow to compute local evaluations explicitly. Before using a symbol defined over a function field $F(X)$, one must verify that its class comes from $\Br(X)$. A class pulled back from $\Br(F)$ has zero sum of local evaluation by reciprocity. Thus, finding a non split algebra over $F$ is not enough to produce a Brauer--Manin obstruction. 

% ----------------------------------------------------------------------
%\section*{Acknowledgments}

% ----------------------------------------------------------------------

\end{document}